\documentclass[a4paper,11pt]{article}
\usepackage[utf8]{inputenc}
\usepackage{fullpage}
\usepackage{xargs}                      
\usepackage[pdftex,dvipsnames]{xcolor}  

\usepackage{fullpage,amsfonts,amsmath,amsthm,amssymb,mathtools,booktabs}
\usepackage{tikz,color,algorithmic,algorithm,multirow,graphicx,caption,subcaption,listings, mathrsfs}
\usepackage[title]{appendix}%
\usepackage{manyfoot}%

\theoremstyle{plain}
\newtheorem{theorem}{Theorem}

\newtheorem{lemma}[theorem]{Lemma}
\newtheorem{remark}[theorem]{Remark}
\newtheorem{example}[theorem]{Example}
\theoremstyle{definition}
\newtheorem{assumption}[theorem]{Assumption}
\newtheorem{definition}[theorem]{Definition}

\newcommand{\R}{\mathbf{R}}

\newcommand{\xs}{x^*}
\newcommand{\kk}{_k}
\newcommand{\kz}{_0}

\newcommand{\kpo}{_{k+1}}

\newcommand{\eqdef}{:=}
\usepackage{soul}

\usetikzlibrary{arrows}

\title{On a Family of Relaxed Gradient Descent Methods for Quadratic Minimization}
\author{Liam MacDonald, Rua Murray and Rachael Tappenden}
\date{}

\begin{document}
\maketitle

\begin{abstract}
    This paper studies the convergence properties of a family of Relaxed $\ell$-Minimal Gradient Descent methods for quadratic optimization; the family includes the omnipresent Steepest Descent method, as well as the Minimal Gradient method. Simple proofs are provided that show, in an appropriately chosen norm, the gradient and the distance of the iterates from optimality converge linearly, for all members of the family. Moreover, the function values decrease linearly, and iteration complexity results are provided. All theoretical results hold when (fixed) relaxation is employed. It is also shown that, given a fixed overhead and storage budget, every Relaxed $\ell$-Minimal Gradient Descent method can be implemented using exactly one matrix vector product. Numerical experiments are presented that illustrate the benefits of relaxation across the family.
\end{abstract}

\paragraph{Keywords.} Steepest descent; relaxation; linear convergence; quadratic optimization; strong convexity; positive definite Hessian.

\section{Introduction}

This work studies a family of gradient descent algorithms for the quadratic minimization problem
\begin{eqnarray}\label{f}
\min_{x\in\R^n} f(x),\quad\text{where}\quad f(x) = \tfrac12 x^TA x - x^Tb,
\end{eqnarray}
where $A\in \R^{n\times n}$ and $x,b \in \R^{n}$. Throughout this work, it is assumed that \eqref{f} is strongly convex so that the unique solution is
\begin{equation}\label{solution}
    \xs = A^{-1}b.
\end{equation}
Problem \eqref{f} is equivalent to solving a system of linear equations  --- the solution to both problems takes the form \eqref{solution} --- so that \eqref{f} is ubiquitous in applied mathematics and the physical and engineering sciences; see, for example \cite{figueiredo2007gradient, huang2008efficient, more1989algorithms, loosli2009quadratic, friedlander1998gradient}. 

Numerous strategies can be employed to solve \eqref{f}, and this work focuses on gradient descent methods. Some of the benefits of gradient-based methods include: they are simple to understand and implement; their iterations are computationally inexpensive; and they are often supported by theoretical convergence results that guarantee a solution to \eqref{f} can be located. (Of course, it is also known that gradient based methods may take many iterations to locate a solution, and they can struggle on ill-conditioned problems.)

In 1847, Cauchy~\cite{cauchy1847methode} presented what is now commonly referred to as the (exact) Steepest Descent (SD) method. Given an initial point $x_0\in \R^n$, for all $k\geq 0$ the iterates of SD take the form
\begin{equation}\label{xupdatenoomega}
    x\kpo = x\kk - \alpha\kk^{\rm SD} g\kk,
\end{equation}
where 
\begin{equation}\label{gradient}
    g\kk = \nabla f(x\kk) = Ax\kk -b,
\end{equation}
is the gradient and 
\begin{equation}\label{steepestStep}
    \alpha\kk^{\rm SD} = \frac{g\kk^Tg\kk}{g\kk^TAg\kk},
\end{equation}
is the step size. Because $f$ is quadratic, $\alpha\kk^{\rm SD}$ is the result of an exact line search, i.e.,
\begin{equation}\label{minproblem}
    \alpha\kk^{\rm SD} = \arg \min\limits_{\alpha} f(x\kk - \alpha g\kk).
\end{equation}

Steepest Descent and its convergence properties have been studied extensively since its inception. The seminal work of Kantorovich \cite{kantorovich1948functional} in the late 1940s establishes his famed inequality\footnote{See \cite{kantorovich1948functional}, bottom of page 142} and uses it to show that for SD the function values converges to the optimum at a linear rate.\footnote{See also \cite{kantorovich1948functional}, the middle of page 144 and the expression $q$ at the bottom of page 145.} Shortly afterwards, in the early 1950s, Akaike~\cite{Akaike1959} characterised the asymptotic behaviour of the method. Akaike showed that, with the choice of stepsize~\eqref{steepestStep}, asymptotically, the search directions alternate within the two-dimensional subspace spanned by the eigenvectors corresponding to the largest and smallest eigenvalues of $A$. This behaviour causes SD to `zig-zag' as it approaches the minimizer, and hence it can be slow on ill-conditioned problems. At a similar time, Forsythe and Motzkin were studying the asymptotic properties of the method, and established several similar results \cite{forsythe1951Asymptoticproperties}. More recently, motivated by the fact that the norm of the gradient is often used as a stopping condition, but that for SD the norm of the gradient can oscillate, Nocedal et al. \cite{Nocedal2002} show that the two-step asymptotic rate of convergence of the norm of the gradient is equal to the one-step asymptotic rate in the function value; this provides understanding of when the norm of the gradient provides a good estimate of the accuracy in the optimal function value.

The quest for step length modifications that reduce or eliminate zig-zagging has been extensive. Current investigations into changing the step size include~\cite{YuanNewStepsize, DaiYuanAnalysis,Huang2021,oviedo2021second}, while alternating between \eqref{alpha} and new `shorter' steps is studied in \cite{Gonzaga}. Other research considers whether it is beneficial to change the step size adaptively, so that the search directions better align with the eigenvectors in later iterations \cite{Zou_2021, Serafino2016}. 
The work \cite{dai2003alternate} consider another `shortened' SD method, 
\cite{forsythe1951Acceleration} investigates a way to `accelerate' SD, and \cite{Raydan2002} develops a new step size --- related to the Barzilai-Borwein \cite{TwoPointBarzilaiBorwein} step length --- which outperforms the Cauchy step \eqref{alpha}. 

Another well studied strategy for reducing zig-zagging is the inclusion of a `relaxation' or `damping' parameter. As early as the 1950s, investigations into the practical behaviour of SD with a (fixed) relaxation parameter had begun \cite{ForsytheIBM}, with a comprehensive set of numerical experiments --- including how the performance changes depending on the choice of the relaxation parameter --- in \cite{stein1951gradient}. Hestenes \cite{hestenes1973solution} discusses `almost-optimum' gradient methods (i.e., relaxed SD), and in \cite{hestenes1990conjugacy} he states that: `We also considered the introduction of a relaxation constant $\beta$ in our algorithm but did not develop an adequate theory for this case'. While the previously mentioned works consider a fixed relaxation parameter, Raydan and Svaiter \cite{Raydan2002} investigate choosing the relaxation parameter randomly at each iteration, and this also shows improvement compared with no relaxation. Importantly, they provide theoretical guarantees, confirming that SD with random relaxation will converge. More recently, Van den Doel and Ascher \cite{Doel2012} show that a fixed relaxation parameter can cause chaotic behaviour.

Closely related to SD is the Minimal Gradient (MG) method \cite{Huang2021,krasnosel1952iteration, DeAsmundis2014efficient}. The iterates of MG also take the form \eqref{xupdatenoomega}, but with $\alpha\kk^{\rm SD}$ replaced by
\begin{equation}\label{MGStep}
    \alpha\kk^{\rm MG} = \frac{g\kk^TAg\kk}{g\kk^TA^2g\kk}.
\end{equation}
The name `Minimal Gradient' method (which appears to be coined in \cite{zhou2006gradient}), comes from the fact that the step size \eqref{MGStep} is that which exactly minimizes the 2-norm of the gradient at the next iterate, i.e., $\alpha\kk^{\rm MG} = \arg \min\limits_{\alpha} \|\nabla f(x\kk - \alpha g\kk)\|_2$. Many of the asymptotic results for SD carry over analogously to the MG method, with \cite{Huang2021} presenting such results rigorously.

In a similar way to the MG method, a class of algorithms exist, where the step sizes arise by exactly
minimizing the gradient in a norm induced by (a power of) the matrix $A$. The goal of this work is to study these algorithms, and provide a comprehensive convergence theory for the entire family. By selecting and working in appropriately chosen norms, the proofs of our theoretical results are straightforward, short, cover all algorithms in the family, and the theory holds even when (fixed) relaxation is included.

\subsection{Notation and Preliminaries}

The following assumption and preliminaries are utilised in this work.

\begin{assumption}\label{Assume1}
The matrix $A\in \R^{n\times n}$ in \eqref{f} is symmetric and positive definite.
\end{assumption}
By Assumption \ref{Assume1}, \eqref{f} is strongly convex. Moreover, because $A$ is positive definite, it has a unique, positive definite square root, denoted by $A^{1/2}$. Throughout this work the convention that $A^0 = I$ is adopted. The eigenvalues of $A$ are denoted by $0<\lambda_n\leq \cdots \leq \lambda_1$, and $\kappa \eqdef \lambda_1/\lambda_n$ is the (2-)\,condition number.

\begin{definition}
Let $B \in \R^{n \times n}$ be a symmetric positive definite matrix. Then define the norm $\|v\|_{B}=\sqrt{v^TBv}$, $v \in \R^n$. 
\end{definition}

\begin{lemma}[Rayleigh Quotient {\cite[Theorem 1.21]{Saad2003}}]\label{lem:Rayleigh}
    For a symmetric matrix $B \in \R^{n \times n}$, and a nonzero vector $v \in \R^n$ it holds that
    \begin{equation}\label{RayleighQuotient}
        \lambda_n(B) \leq \frac{v^TBv}{v^Tv} \leq \lambda_1(B).
    \end{equation}
\end{lemma}

\begin{theorem}[Kantorovich Inequality~{\cite[Lemma~5.8]{Saad2003}}]\label{Thm:Kantorovich}
Let $B$ be any symmetric positive definite matrix and $\lambda_1$ and $\lambda_n$ its largest and smallest eigenvalues, respectively. Then
\begin{eqnarray}\label{Kantorovich}
\frac{(x^TBx) (x^TB^{-1}x)}{\|x\|_2^4} \leq \frac{(\lambda_1 + \lambda_n)^2}{4\lambda_1 \lambda_n},\quad \forall \; x \neq 0.
\end{eqnarray}\
\end{theorem}

\subsection{Outline}

The remainder of this work is organized as follows. In Section~\ref{sec:FamilyMethods}, the class of algorithms to be studied is defined, and the contributions of this work are stated. Section~\ref{sec:ConvergenceProofs} presents convergence theory for the algorithms under consideration, including convergence results for the norm (induced by a power of $A$) of the gradient, for the distance of the iterates from optimality, function values, and iteration complexity results are also provided. In Section~\ref{sec:EfficientComp} it is shown that, given a fixed overhead and memory budget, every algorithm in the family can be implemented using one matrix vector product per iteration. Finally, numerical experiments are presented in Section~\ref{sec:NumericalResults}.

\section{A family of relaxed gradient descent methods}\label{sec:FamilyMethods}

In this section the class of algorithms considered throughout this work is described. Before stating the definition, note that straightforward algebraic manipulation yields
\begin{align}\label{eq:fvalsvsg}
    f(x) - f^* = \tfrac{1}{2}(x-x^*)^TA(x-x^*) = \tfrac{1}{2}\|g(x)\|_{A^{-1}}^2.
\end{align}
Recall that for SD, the iterates take the form \eqref{xupdatenoomega}, where the step size \eqref{steepestStep} is found via an exact line search \eqref{minproblem}. Thus, employing the relationship \eqref{eq:fvalsvsg} shows that \eqref{steepestStep} is equivalently found via $\alpha\kk^{\rm SD} = \arg \min\limits_{\alpha} \|\nabla f(x\kk - \alpha g\kk)\|_{A^{-1}}$. Moreover, it has already been seen that the step length for the MG method is computed similarly, but w.r.t. the 2-norm. This motivates the following lemma.

\begin{lemma}\label{lem:alpha}
Let $f$ be defined in \eqref{f}, let Assumption \ref{Assume1} hold, and let $\ell \in \{0,\tfrac12,1,\tfrac32,2,\tfrac52,3,\dots\}$ be fixed. Given $x\kk\in \R^n$, let $g\kk$ be given in \eqref{gradient} and define
\begin{eqnarray}\label{alpha}
\alpha\kk = \frac{y\kk^Ty\kk}{y\kk^TAy\kk}, 
\end{eqnarray}
where
\begin{eqnarray}\label{y}
y\kk = A^{\ell} g\kk.
\end{eqnarray}
Then $\alpha\kk = \arg \min\limits_{\alpha}\|\nabla f(x_k-\alpha\ g_k)\|_{A^{2\ell-1}}$.
\end{lemma}
\begin{proof}
Observe that
\begin{eqnarray}
\|\nabla f(x_k-\alpha\ g_k)\|_{A^{2\ell-1}}^2 
 &=& (g\kk - \alpha A g\kk)^TA^{2\ell-1}(g\kk - \alpha A g\kk)\notag\\
&=& g\kk^T A^{2\ell-1}g\kk - 2\alpha g\kk^T A^{2\ell}g\kk  + \alpha^2 g\kk^TA^{2\ell+1}g\kk\notag\\
&=& y\kk^T A^{-1}y\kk - 2\alpha y\kk^T y\kk  + \alpha^2 y\kk^TAy\kk\label{gkpoinA2lm1}.
\end{eqnarray}
Differentiating w.r.t. $\alpha$, and rearranging, gives \eqref{alpha}.
\end{proof}
Notice that using $\ell = 0$ in Lemma~\ref{lem:alpha} recovers the step size for SD \eqref{steepestStep}, while using $\ell = 1/2$ (and recalling the convention that $A^0 = I$) recovers the step size for the MG method \eqref{MGStep}.

\begin{remark}
    It is important to stress here that the fractional matrix powers appearing in \eqref{y} are purely for notational convenience. A matrix square root need never be computed in practice; the calculations can always be arranged to avoid this (see also Section~\ref{sec:EfficientComp}).
\end{remark}

As previously mentioned, relaxation can be beneficial in terms of the practical performance of an algorithm. Thus, consider an iterative process of the form $x\kpo = \Phi (x\kk)$, where $\Phi : \R^n \to \R^n$. Let $\omega \in (0,2)$, and consider the update defined via $x\kpo = (1-\omega)x\kk + \omega \Phi (x\kk).  $
Applying this `acceleration' process to the update in \eqref{xupdatenoomega} gives
\begin{eqnarray*}\label{relaxedGDiterates2}
x\kpo = (1-\omega)x\kk + \omega (x\kk - \alpha\kk g\kk) = x\kk - \omega \alpha\kk g\kk.
\end{eqnarray*}

Algorithms whose iterates take the form above are often referred to as relaxed (or damped) gradient descent methods, where $\omega$ is the relaxation/damping parameter. The class of algorithms that is studied in this work is defined now. 

\begin{definition}[Relaxed $\ell$-Minimal Gradient Descent (Relaxed $\ell$-MGD) Methods]\label{def:Algorithms}
   Let $f$ be given in \eqref{f}, let Assumption \ref{Assume1} hold, let $\ell \in \{0,\tfrac12,1,\tfrac32,2,\tfrac52,3,\dots\}$ be fixed, and choose $\omega \in (0,2)$. Given an initial point $x\kz \in \R^n$, the $k$-th iteration of a Relaxed $\ell$-Minimal Gradient Descent Method is defined by 
   \begin{equation}\label{xupdate}
       x\kpo = x\kk - \omega \alpha_k g\kk,
   \end{equation}
   where $g\kk$ and $\alpha\kk$ are given in \eqref{gradient} and \eqref{alpha}, respectively.
\end{definition}

Definition~\ref{def:Algorithms} is adapted from Definition~1 in \cite{Pronzato_2005}, which considers a family of methods that they call $P$-gradient Algorithms.\footnote{Definition~1 in \cite{Pronzato_2005} is stated in full in Appendix~\ref{app:PronzatoEquiv} for ease of reference.} If $\omega = 1$ (i.e., no relaxation), then the (Relaxed) $\ell$-MGD methods described in Definition~\ref{def:Algorithms} belong to the class of $P$-gradient Algorithms. However, \eqref{xupdate} includes a relaxation parameter $\omega \in (0,2)$, so that any algorithm satisfying Definition~\ref{def:Algorithms} with $\omega \neq 1$ is not a $P$-gradient Algorithm. It can be seen that (relaxed) SD ($\ell = 0$) and the (relaxed) MG method ($\ell = 1/2$) are Relaxed $\ell$-MGD methods (and if $\omega=1$ then they are also $P$-gradient Algorithms).

Pronzato et al. \cite{Pronzato_2005} analyze the $P$-gradient Algorithms over a Hilbert space and provide a comprehensive asymptotic convergence theory for the whole class. In essence, they generalize the work of Akaike~\cite{Akaike1959} for SD, to the family of $P$-gradient Algorithms. The work \cite{Huang2021} also studies $P$-gradient algorithms, although their convergence theory holds over $\R^n$.

To the best of our knowledge, a thorough non-asymptotic theoretical study of the family of relaxed $\ell$-MGD methods (Definition~\ref{def:Algorithms}) is lacking. A key contribution of this work is to fill this gap.

Before continuing, it is helpful to understand why relaxation can benefit the practical performance of (relaxed) $\ell$-MGD methods. So, suppose that $g\kk \neq 0$ (i.e., suppose the current iterate $x\kk$ is not the solution to \eqref{f}). Then
\begin{eqnarray}\label{zigzag}
g\kpo^TA^{2\ell}g\kk &=& g\kk^T(I-\omega \alpha\kk A)A^{2\ell}g\kk\notag\\
 &=& g\kk^TA^{2\ell}g\kk- \omega \alpha\kk g\kk^T A^{2\ell+1} g\kk\notag\\
 &=& g\kk^TA^{2\ell}g\kk- \omega\frac{g\kk^T A^{2\ell}g\kk}{g\kk^TA^{2\ell+1}g\kk } g\kk^T A^{2\ell+1} g\kk\notag\\
 &=& (1-\omega)\|g\kk\|_{A^{2\ell}}^2.
\end{eqnarray}
Moreover, combining \eqref{y} with \eqref{zigzag} shows that 
\begin{equation}\label{zigzagy}
    y\kpo^Ty\kk = (1-\omega)\|y\kk\|_{2}^2.
\end{equation}

In words, if $\omega = 1$, then \eqref{zigzag} shows that consecutive gradients $g\kk$ and $g\kpo$, generated by any algorithm satisfying Definition~\ref{def:Algorithms} are $A^{2\ell}$-conjugate, while \eqref{zigzagy} shows that the vectors $y\kk$ and $y\kpo$ are orthogonal. Thus, \eqref{zigzag} and \eqref{zigzagy} confirm what is known in practice, that algorithms belonging to this class (without relaxation), tend to zig-zag near the solution, so their practical behaviour can be poor, especially for ill-conditioned problems with highly elliptical contours. Importantly, \eqref{zigzag} and \eqref{zigzagy} also show that if $\omega \neq 1$, then consecutive gradients, $g\kk$ and $g\kpo$ are no longer $A^{2\ell}$-conjugate (and consecutive $y\kk$s are not orthogonal), so that relaxation helps to `break' the zig-zagging and push the iterates toward the minimizer.

Finally, notice that for any vector $y\kk\neq 0$, $\alpha\kk$ in \eqref{alpha} takes the form of the reciprocal of the Rayleigh quotient so that the step length is an approximation to the reciprocal of an eigenvalue of $A$. Thus $\alpha\kk \in [1/\lambda_1,1/\lambda_n]$ and $\alpha\kk>0$ by Assumption~\ref{Assume1}.

\subsection{Contributions}
The main contributions of this work are summarized now (in no particular order).
\begin{enumerate}
    \item \emph{Simplified convergence proofs.} By working with norms induced by an appropriate power of $A$, short and simple convergence proofs are provided for all relaxed $\ell$-MGD methods. Specifically, the norm of the gradient, and the distance of the iterates from optimality, converge linearly.
    \item \emph{Linear convergence of function values.} To the best of our knowledge, this is the first work to establish linear convergence of the function values for the family of relaxed $\ell$-MGD methods. Moreover, a counterexample is provided to confirm that all algorithms from Definition~\ref{def:Algorithms} with $\ell \geq 1/2$ can exceed the known rate $(\kappa-1)^2/(\kappa+1)^2$ of SD.
    \item \emph{Iteration complexity.} The results in this work are non-asymptotic. Thus, we are able to provide iteration complexity results for the relaxed $\ell$-MGD methods, i.e., we provide an explicit expression for the number of iterations $K$ needed to ensure convergence to a given stopping tolerance.
    \item \emph{Relaxation.} All theoretical results hold when fixed relaxation $\omega \in (0,2)$ is employed.
    \item \emph{Computational cost.} Given an initial overhead and storage capacity of $\lfloor\ell\rfloor + 2$ matrix vector products, respectively, every algorithm satisfying Definition~\ref{def:Algorithms} uses exactly one matrix vector product per iteration.
\end{enumerate}

\section{Convergence Properties}\label{sec:ConvergenceProofs}
In this section, theoretical convergence guarantees for the Relaxed $\ell$-MGD methods (Definition~\ref{def:Algorithms}) are consolidated. A key component of this work is that by studying quantities in an appropriate norm (induced by a power of $A$) results can be stated in a way that is transparent, that allows for ease of comparison between different methods in the class, and the proofs are concise (the Cauchy-Schwarz inequality is the main workhorse). Importantly, the results here are both non-asymptotic and they are applicable when fixed relaxation is used, filling a gap in the current literature.

\subsection{Convergence in the $A^{2\ell-1}$-norm}

This study begins with a result showing that the gradient and distance of the iterates from optimality both converge linearly when measured with respect to an appropriately chosen norm. This holds for all algorithms in the class, and includes fixed relaxation. In the absence of relaxation, results similar to those stated below can be found in the literature \cite{Pronzato_2005}, although careful reading may be required to recognise them. To the best of our knowledge, a formal statement of the results, which includes relaxation, and a proof in the setting of $\R^n$ and with the `power of $A$' norms does not appear elsewhere in the literature.

\begin{theorem}
\label{Thm:gmonodecreaserelaxed}
Let $f$ be given in \eqref{f}, let Assumption \ref{Assume1} hold, and fix $\ell \in \{0,\tfrac12,1,\tfrac32,2,\tfrac52,3,\dots\}$ and $\omega \in (0,2)$. Given an initial point $x_0\in \R^n$, for $k\geq 0$, let the iterates be given in Definition~\ref{def:Algorithms} and define
\begin{eqnarray}\label{gmonodecreaserelaxedc}
c(\omega) \eqdef \left(1 - \omega(2-\omega)\frac{4\lambda_1 \lambda_n}{(\lambda_1+\lambda_n)^2}\right).
\end{eqnarray}
Then, 
\begin{eqnarray}\label{gmonodecreaserelaxed}
\|g\kk\|_{A^{2\ell-1}}^2 \leq (c(\omega))^k \, \|g\kz\|_{A^{2\ell-1}}^2,
\end{eqnarray}
\begin{eqnarray}\label{xvxs}
\|x\kk-\xs\|_{A^{2\ell +1}}^2 \leq (c(\omega))^k \; \|x\kz-\xs\|_{A^{2\ell +1}}^2,
\end{eqnarray}
and 
\begin{equation}\label{yconverge}
    \|y\kk\|_{A^{-1}}^2 \leq (c(\omega))^k \;\|y\kz\|_{A^{-1}}^2,
\end{equation}
\end{theorem}
\begin{proof}
Consider \eqref{gmonodecreaserelaxed}, and observe that for all $k\geq 0$, $\phi_k(\omega) \eqdef \|g\kpo\|_{A^{2\ell-1}}^2  = \|g\kk - \omega \alpha\kk Ag\kk\|_{A^{2\ell-1}}^2$ is quadratic in $\omega$, obtaining its global minimizer at $\omega = 1$. Therefore, by symmetry, $\phi_k(0) = \phi_k (2)$, and for any $\omega \in (0,2)$, it holds that $\phi_k(\omega) < \phi_k(0)$. 

Now, 
\begin{eqnarray}
\|g\kpo\|_{A^{2\ell-1}}^2 
&\overset{\eqref{gkpoinA2lm1}}{=}& 
\|g\kk\|_{A^{2\ell-1}}^2 - 2\omega \alpha\kk \|y\kk\|_2^2 + \omega^2 \alpha\kk^2\|y\kk\|_A^2\notag\\
&\overset{\eqref{alpha}}{=}& \|g\kk\|_{A^{2\ell-1}}^2 - 2\omega\frac{\|y\kk\|_2^2}{\|y\kk\|_A^2} \|y\kk\|_2^2 + \omega^2\frac{\|y\kk\|_2^4}{\|y\kk\|_A^4}\|y\kk\|_A^2\notag\\
&=& \|g\kk\|_{A^{2\ell-1}}^2 - \omega(2-\omega)\frac{\|y\kk\|_2^4}{\|y\kk\|_A^2}\notag\\
&=& \left(1 - \omega(2-\omega)\frac{\|y\kk\|_2^4}{\|y\kk\|_A^2 \|g\kk\|_{A^{2\ell-1}}^2}\right)\|g\kk\|_{A^{2\ell-1}}^2\notag\\
&\overset{\eqref{y}}{=}& \left(1 - \omega(2-\omega)\frac{\|y\kk\|_2^4}{\|y\kk\|_A^2 \|y\kk\|_{A^{-1}}^2}\right)\|g\kk\|_{A^{2\ell-1}}^2 \label{PronzatoEquiv}\\
&\overset{\eqref{Kantorovich}}{\leq}& \left(1 - \omega(2-\omega)\frac{4\lambda_1 \lambda_n}{(\lambda_1+\lambda_n)^2}\right)\|g\kk\|_{A^{2\ell-1}}^2.\label{relaxedgraddecrease}
\end{eqnarray}
Note that $4\lambda_1 \lambda_n/(\lambda_1 + \lambda_n)^2\leq1$ (with equality iff $\lambda_1=\lambda_n$), and for all $\omega \in (0,2)$ it holds that $\omega (2-\omega)\leq 1$ (with equality iff $\omega =1$), so $c(\omega) \leq 1$. Repeated application of \eqref{relaxedgraddecrease} gives \eqref{gmonodecreaserelaxed}.

Finally, \eqref{xvxs} follows from
\begin{eqnarray}\label{xdifftog}
\|x\kk-\xs\|_{A^{2\ell +1}}^2 = \|A(x\kk-\xs)\|_{A^{2\ell -1}}^2 = \|g\kk\|_{A^{2\ell -1}}^2,
\end{eqnarray}
while \eqref{yconverge} follows from $\|g\kk\|_{A^{2\ell - 1}}^2 = g\kk^TA^{2\ell - 1}g\kk \overset{\eqref{y}}{=} y\kk^T A^{-1}y\kk = \|y\kk\|_{A^{-1}}^2$.
\end{proof}

Theorem~\ref{Thm:gmonodecreaserelaxed} shows that the gradient converges linearly (w.r.t. the  $\|\cdot\|_{A^{2\ell-1}}$-norm) for all Relaxed $\ell$-MGD methods. It is known that for SD ($\ell = 0$ and $\omega=1$), the 2-norm of the gradient can oscillate. However, Theorem~\ref{Thm:gmonodecreaserelaxed} confirms that for SD, the gradient decreases monotonically w.r.t. the $A^{-1}$-norm.

Note that when $\omega = 1$ (i.e., no relaxation), \eqref{gmonodecreaserelaxedc} becomes the known rate
\begin{equation}\label{ratenorelax}
    c(1) = \left(1 - \frac{4\lambda_1 \lambda_n}{(\lambda_1+\lambda_n)^2}\right) = \left(\frac{\lambda_1 - \lambda_n}{\lambda_1 + \lambda_n}\right)^2.
\end{equation}

Moreover, the theoretical rate $c(1)$ in \eqref{ratenorelax} is better than if $\omega \neq 1$ (no relaxation). Nevertheless, relaxation often improves the practical behaviour of gradient descent methods (see Section~\ref{sec:NumericalResults}). Finally, note that the rate given in \eqref{ratenorelax} matches the rate of convergence given by Pronzato et al. \cite{Pronzato_2005} (see Appendix~\ref{app:PronzatoEquiv} for further details).

\subsection{Convergence of the function values}
\label{sec:funcvals}

A major contribution of this work is to establish convergence of the function values for the family of algorithms in Definition~\ref{def:Algorithms}. It is known that for SD (without relaxation) the function values evolve at the rate $c(1)$. However, this section begins with a counterexample, which shows that the (relaxed) $\ell$-MGD methods with $\ell \geq 1/2$ can exceed the rate \eqref{ratenorelax}.

Following the counterexample, two different but complementary (linear) rates of decay for the Relaxed $\ell$-MGD methods are given. The first shows that the function values decay in the same way for all algorithms, which supports numerical experiments that show the algorithms perform similarly in practice (see also Section~\ref{sec:NumericalResults}). The second involves the familiar `Kantorovich' rate, but there is also a constant factor (a power of the condition number), and this leads to a better iteration complexity bound (see also Section~\ref{sec:IterationComplexity}).

\subsubsection{A counterexample}\label{sec:Counterexample}

Consider Exercise 2.11 in Fletcher \cite{fletcher2000practical}, where the given function $f(x) = 10x_1^2 + x_2^2$ can be written in the form \eqref{f} with $b = 0$ and $A = {\rm diag}(\lambda_1,\lambda_2)$, and where the eigenpairs are $(\lambda_1,v_1)=(20,e_1)$ and $(\lambda_2,v_2)=(2,e_2)$. The initial point is given as $x\kz = (1/10, 1)^T$, and it is this point that forces the worst case behaviour in SD. In general, for any iteration $k\geq 0$ of a (relaxed) $\ell$-MGD method, the gradient can be expressed in the basis of eigenvectors as $g\kk = \sum_{i=1}^n = c_i v_i$, for some constants $c_i$, $\forall i$. Here, $g\kz = Ax\kz = 2 e_1 + 2 e_2$, so that the initial point $x\kz$ forces the weighting apportioned to each of the eigenvectors to be equal, i.e., $c_1 = c_2 = 2$. 

\begin{figure}[h]
    \centering
    \includegraphics[width=0.5\textwidth]{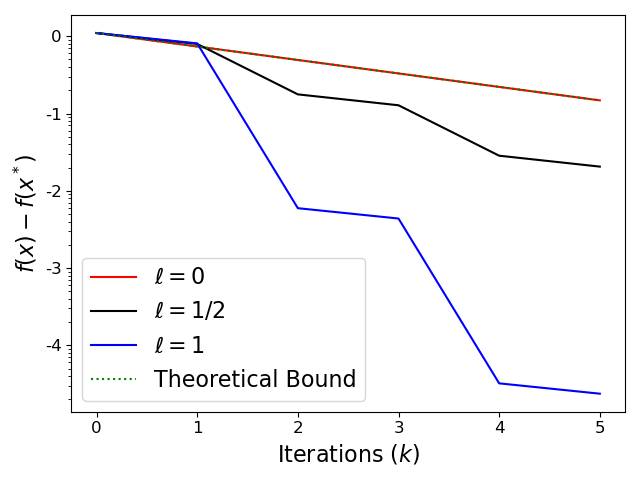}
    \caption{Function values for first 5 iterates of a (relaxed) $\ell$-MGD method applied to $f(x) = 10x_1^2 + x_2^2$, from an initial point $x\kz = (1/10,1)^T$ with $f(x_0) = 1.1$.}
    \label{fig:CounterExample}
\end{figure}

Now consider Table~\ref{tab:CounterExample} and Figure~\ref{fig:CounterExample} and note that $f(x\kz)=1.1$. Table~\ref{tab:CounterExample} shows that the function values $f(x_1)$ for the $\ell=1/2$ and $\ell = 1$ (relaxed) $\ell$-MGD methods are larger than the known theoretical bound for SD (while SD meets its bound). Figure~\ref{fig:CounterExample} shows that SD maintains the worst case rate for this example for all iterations. However, Figure~\ref{fig:CounterExample} also shows that, apart from the first iteration, the function values for the $\ell=1/2$ and $\ell = 1$ (relaxed) $\ell$-MGD methods are lower than for SD.

Indeed, all (relaxed) $\ell$-MGD methods will perform worse (in terms of the function value) than SD in the first iteration. Hence the `Kantorovich' rate \eqref{ratenorelax} cannot be guaranteed for the \emph{one-step} function values.

\begin{table}[h]
    \centering
    \begin{tabular}{||c|c|c|c|c||}
    \hline
        $k$ & $\ell=0$ &$\ell=1/2$& $\ell=1$&Theoretical Bound \\
        \hline
         $f(x_{1})$& 0.7364& 0.7948& 0.8084&0.7634 \\
         $f(x_2)$& 0.4929&0.1769&0.0060&0.4929\\
         \hline
    \end{tabular}
    \caption{Function values for first 2 iterates of a (relaxed) $\ell$-MGD method applied to $f(x) = 10x_1^2 + x_2^2$, from an initial point $x\kz = (1/10,1)^T$ with $f(x_0) = 1.1$.}
    \label{tab:CounterExample}
\end{table}

\begin{remark}
    In general, for any Relaxed $\ell$-MGD method, the worst case rate of $\|g\kk\|_{A^{2\ell-1}}^2$ is obtained by considering $y\kk$ expressed in a basis of eigenvectors, and choosing an initial point $x\kz\in \R^n$ that forces the coefficients of the $v_1$ and $v_n$ terms to be equal. For example, here, a `worst case' initial point for the MG method is $x\kz = (\frac1{10\sqrt{10}},1)^T$, while a `worst case' initial point for the Relaxed $1$-MGD method is $x\kz = (1/100,1)^T$.
\end{remark}

\subsubsection{Linear Convergence of the function values}
\label{Sec:functionvals}

We are now ready to state one of the major results of this work, which shows that function values decay at a linear rate, for all Relaxed $\ell$-Minimal Gradient algorithms described in Definition~\ref{def:Algorithms}. The proof relies on several technical results that are given in Lemmas~\ref{Lemmainduction} and \ref{lemmaBounding} in Appendix~\ref{sec:technicallemmas}.

\begin{theorem}\label{thm:fvaldecrease}
Let $f$ be given in \eqref{f}, let Assumption~\ref{Assume1} hold, fix $\ell \in \{0,\tfrac12,1,\tfrac32,2,\tfrac52,3,\dots\}$ and $\omega \in (0,2)$, and let $c(\omega)$ be given in \eqref{gmonodecreaserelaxedc}. Given an initial point $x_0\in\R^n$, for $k\geq 0$, let the iterates be given in Definition~\ref{def:Algorithms}. Then, 
\begin{enumerate}
    \item Rate I: \begin{equation}\label{eq:fvaldecrease}
     f(x\kk)-f^* \leq \left(1- \omega(2-\omega) \frac{1}{\kappa}\right)^k(f(x\kz)-f^*);
\end{equation}
\item Rate II: 
\begin{eqnarray}\label{fdecreaseLine}
    f(x_k)-f^* \leq 
    \kappa^{2\ell}(c(\omega))^{k} \, (f(x_0)-f^*).
    \end{eqnarray}
\end{enumerate}
\end{theorem}
\begin{proof}
We begin by establishing Rate I. Similarly to the proof of Theorem~\ref{Thm:gmonodecreaserelaxed}, and recalling~\eqref{eq:fvalsvsg}, observe that for any $k\geq 0$, $\varphi_k(\omega) \eqdef \|g\kpo\|_{A^{-1}}^2  = \|g\kk - \omega \alpha\kk Ag\kk\|_{A^{-1}}^2$ is a quadratic in $\omega$, which obtains its global minimizer at $\omega = \|g\kk\|_2^2/(\alpha\kk \|g\kk\|_A^2) \geq 1$. Therefore, and for all $\omega \in (0,2)$, it holds that $\phi_k(\omega) < \phi_k(0)$. Now, to prove \eqref{eq:fvaldecrease} observe that
    \begin{eqnarray}
\|g\kpo\|_{A^{-1}}^2 
&=& \|g\kk - \omega \alpha\kk Ag\kk \|_{A^{-1}}^2 \notag\\
&=& \|g\kk\|_{A^{-1}}^2  -2\omega \alpha\kk \|g\kk\|^2 + \omega^2 \alpha\kk^2 \|g\kk\|_{A}^2  \notag\\
&=& \|g\kk\|_{A^{-1}}^2  -2\omega \frac{\|y\kk\|_2^2}{\|y\kk\|_A^2} \|g\kk\|^2 + \omega^2 \frac{\|y\kk\|_2^4}{\|y\kk\|_A^4} \frac{\|g\kk\|_2^2}{\|g\kk\|_2^2}\|g\kk\|_{A}^2  \notag\\
&\overset{\eqref{intermediateThmresult}}{\leq}& \|g\kk\|_{A^{-1}}^2  -2\omega \frac{\|y\kk\|_2^2}{\|y\kk\|_A^2} \|g\kk\|^2 + \omega^2 \frac{\|y\kk\|_2^2}{\|y\kk\|_A^2 } \|g\kk\|_2^2  \notag\\
&\leq& \left(1  -\omega(2-\omega) \frac{\|y\kk\|_2^2}{\|y\kk\|_A^2} \frac{\|g\kk\|^2}{\|g\kk\|_{A^{-1}}^2}\right) \|g\kk\|_{A^{-1}}^2.\notag\\
&\overset{\eqref{lem:Rayleigh}}{\leq}& \left(1  -\omega(2-\omega) \frac{\lambda_n}{\lambda_1}\right) \|g\kk\|_{A^{-1}}^2.\notag
\end{eqnarray}

To establish Rate II, it holds that
\begin{eqnarray*}
        f(x\kk) - f^* 
        &\overset{\eqref{eq:fvalsvsg}}{=}& \tfrac{1}{2}(x\kk-x^*)^TA(x\kk-x^*)\\
        &\overset{\eqref{ineq1}}{\leq}& \frac12\frac{1}{\lambda_n^{2\ell}}  \|x\kk-x^*\|_{A^{2\ell+1}}^2\notag \\
        &\overset{\eqref{xdifftog}}{=}& \frac{1}{2}\frac{1}{\lambda_n^{2\ell}} \; (c(\omega))^{k}  \; \|x\kz-x^*\|_{A^{2\ell+1}}^2\notag\\
        &\overset{\eqref{ineq2}}{\leq}& \left(\frac{\lambda_1}{\lambda_n}\right)^{2\ell} \; (c(\omega))^{k}  \; (f(x\kz)-f^*).
    \end{eqnarray*}
    The proof of \eqref{fdecreaseLine} is complete.
\end{proof}

Theorem~\ref{thm:fvaldecrease} shows that for all Relaxed $\ell$-MG methods, the function values decrease linearly. Rates I and II are in a sense, complementary. Notice that \eqref{eq:fvaldecrease} is independent of $\ell$, so that Rate~I suggests that evolution of the function values is the same for every Relaxed $\ell$-MG method. This appears to be mirrored in the numerical experiments, which show that the algorithms perform similarly in practice, regardless of $\ell$; see also Section~\ref{sec:NumericalResults}. At the same time, the rate in \eqref{eq:fvaldecrease} is pessimistic, and it is known that SD achieves the better rate given in \eqref{gmonodecreaserelaxedc}.

On the other hand, the rate in \eqref{fdecreaseLine} is more optimistic, matching that in \eqref{gmonodecreaserelaxedc}, and when $\ell=0$, it recovers the known rate for SD. However, Rate II involves a constant that is a power of the condition number, and this constant is large when $\ell$ is large. Moreover, because Rate II depends on $\ell$, it also suggests that there is a difference in performance when comparing the Relaxed $\ell$-MG methods (i.e., it suggests worse performance for algorithms with larger $\ell$).

\subsection{Convergence in the 2-norm}

Typically, it is desirable to develop convergence theory for algorithms that measures the 2-norm of the gradient, because this is what is usually computed in practice. This section adapts the results from Theorem~\ref{thm:fvaldecrease} to results depending on $\|g\kk\|_2$.

\begin{lemma}\label{lem:g2normdecrease}
    Let the conditions of Theorem~\ref{thm:fvaldecrease} hold. Then, 
    \begin{enumerate}
        \item Rate I: 
        \begin{equation}\label{eq:g2normdecrease}
        \|g\kk\|_2^2 \leq \kappa \left(1- \omega(2-\omega) \tfrac{1}{\kappa}\right)^k\|g\kz\|_2^2,
    \end{equation}
    \item Rate II:
    \begin{enumerate}
    \item if $\ell = 0$, then 
    \begin{eqnarray}\label{gmonodecreaserelaxed2normell0}
\|g\kk\|_{2}^2 \leq \kappa \, (c(\omega))^k \, \|g\kz\|_{2}^2;
\end{eqnarray}
\item and if $\ell \in \{\frac12,1,\tfrac32,2,\tfrac52,3,\dots\}$, then 
\begin{eqnarray}\label{gmonodecreaserelaxed2norm}
\|g\kk\|_{2}^2 \leq \kappa^{2\ell-1} (c(\omega))^k \, \|g\kz\|_{2}^2.
\end{eqnarray}
\end{enumerate}
    \end{enumerate}   
\end{lemma}
\begin{proof}
    To show Rate I, by strong convexity (Assumption~\ref{Assume1}), the function values can be sandwiched as follows,
    $\frac{1}{2\lambda_1} \|g\kk \|^2_2 \leq f(x_k) - f(x^*) \leq \frac{1}{2\lambda_n} \|g\kk\|_2^2$. 
Combining this with \eqref{eq:fvaldecrease} gives \eqref{eq:g2normdecrease}.\\

Now, to establish Rate II(a), \eqref{RayleighQuotient} shows that $(1/\lambda_1)\|g\kk\|_2^2 \leq \|g\kk\|_{A^{-1}}^2$ and $\|g\kz\|_{A^{-1}}^2 \leq (1/\lambda_n)\|g\kz\|_2^2$. Combining this with \eqref{gmonodecreaserelaxed}, and rearranging, gives the result \eqref{gmonodecreaserelaxed2normell0}.

    For Rate II(b), \eqref{RayleighQuotient} also shows that, $\lambda_n^{2\ell-1} \|g\kk\|_{2}^2 \leq \|g\kk\|_{A^{2\ell-1}}^2 $ and $\|g\kz\|_{A^{2\ell-1}}^2 \leq \lambda_1^{2\ell-1} \|g\kz\|_{2}^2$. Combining this with \eqref{gmonodecreaserelaxed}, and rearranging, gives the result \eqref{gmonodecreaserelaxed2norm}.
\end{proof}

\begin{remark}
    Note that $\|g\kk\|_2^2 = \|x\kk-\xs\|_{A^2}^2 $, so that \eqref{eq:g2normdecrease}, \eqref{gmonodecreaserelaxed2normell0} and \eqref{gmonodecreaserelaxed2norm} can be directly translated into results measuring the distance of the iterates from optimality in the $A^2$-norm.
\end{remark}

It is known that Steepest Descent (Definition~\ref{def:Algorithms} with $\ell=0$ and $\omega = 1$)  exhibits oscillations in the 2-norm of the gradient. In particular, \cite[Theorem 4.1]{Nocedal2002} shows that

    \begin{equation}\label{nocedalgoscillate}
        \frac{\|g\kpo\|_2^2}{\|g\kk\|_2^2} \leq \frac{(\kappa-1)^2}{4\kappa},
    \end{equation}
     so that if $\kappa > 3 + 2\sqrt{2}$, then oscillations can occur. Importantly, \eqref{nocedalgoscillate} and Lemma~\ref{lem:g2normdecrease} are compatible. Notice that as $\kappa \to \infty$, \eqref{nocedalgoscillate} tends to $\kappa/4$, i.e., from one iteration to the next, the 2-norm of the gradient may increase by $\approx \kappa/4$. The bound \eqref{gmonodecreaserelaxed2normell0} is not an `iteration-to-iteration' rate. Instead, at every iteration the \emph{upper bound} $\kappa \|g\kz\|_2^2$ (note the inclusion of the condition number $\kappa$), is pushed down at the rate $c(1)$ (assuming no relaxation), which eventually drives the 2-norm of the gradient to zero. It is also this constant factor $\kappa$ that allows for the iteration-to-iteration oscillations in the 2-norm to occur (which is consistent with that described in \cite{Nocedal2002}).

\subsection{Iteration complexity results}\label{sec:IterationComplexity}

Finally, iteration complexity results for the family of relaxed $\ell$-MGD methods are presented. These provide an explicit expression for $K$, the number of iterations required to push $f(x\kk) - f^*$ below some desired tolerance $\epsilon>0$.

\begin{theorem}\label{Thm:ComplexityOne}
Let the conditions of Theorem~\ref{thm:fvaldecrease} hold and let $\epsilon,\hat{\epsilon} >0$. Then, for the Relaxed $\ell$-MG methods (Definition~\ref{def:Algorithms}),
\begin{enumerate}
    \item $f(x_{K})-f^* \leq \epsilon$, where 
    \begin{equation}\label{eq:K}
        K > \frac{\kappa}{\omega(2-\omega)}\left(\ln{\frac{f(x\kz)-f^*}{\epsilon}} \right),
    \end{equation}
\item or, $f(x_{\hat K})-f^* \leq \hat{\epsilon}$, where 
\begin{equation}\label{eq:Khat}
    \hat{K} > \frac{1}{\omega(2-\omega)}\frac{(\kappa+1)^2}{4\kappa}\left(\ln{\frac{\kappa^{2\ell}\left(f(x\kz)-f^*\right)}{\hat{\epsilon}}} \right).
\end{equation}
\end{enumerate}
\end{theorem}
\begin{proof}
    Note that $(1-c)^{1/c} \leq e^{-1}$. Then, to establish (1) consider
    \begin{eqnarray*}
        f(x_{K})-f^* &\overset{\eqref{eq:fvaldecrease}}{\leq}& \left(1-  \frac{\omega(2-\omega)}{\kappa}\right)^K(f(x\kz)-f^*)\\
        &<&  \left(\left(1-  \frac{\omega(2-\omega)}{\kappa}\right)^{\frac{\kappa}{\omega(2-\omega)} }\right)^{  \left(\ln{\frac{f(x\kz)-f^*}{\epsilon}}\right) } (f(x\kz)-f^*)\\
        &\leq& (e^{-1})^{\left(\ln{\frac{f(x\kz)-f^*}{\epsilon}}\right)} (f(x\kz)-f^*)\\
         &= & \epsilon.
    \end{eqnarray*}
    The proof of (2) is similar, so is omitted.
\end{proof}

\begin{remark}
    In an analogous way, one can obtain iteration complexity bounds involving, for example, (i) $\|g_K\|_2^2$ (for both Rate I and Rate II); (ii) $\|x_K-\xs\|_{A^2}^2 $ ; (iii) $\|g_K\|_{A^{2\ell-1}}^2$; (iv) $\|x_K-\xs\|_{A^{2\ell+1}}^2$; or (v) $\|y_K\|_{A^{-1}}^2$. The results and proofs are essentially identical to those in Theorem~\ref{Thm:ComplexityOne}, so are omitted for brevity.
\end{remark}

Theorem~\ref{Thm:ComplexityOne} shows that $O(\ln \frac{1}{\epsilon})$ iterations are needed for convergence (and this holds true for \emph{any} Relaxed $\ell$-MGD method). Notice that the iteration complexity bound is \emph{worse} if relaxation is used (i.e., when $\omega \neq 1)$ although this is usually not the case in practice (as will be shown in the numerical experiments in Section~\ref{sec:NumericalResults}).

\begin{figure}
    \centering
    \includegraphics[height=0.25\textheight]{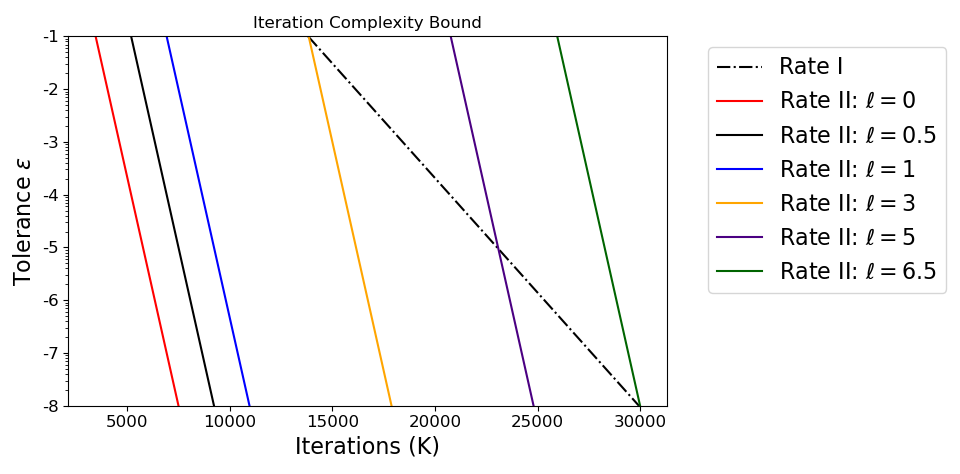}
    \caption{Plot illustrating how the iteration complexity bound in Theorem~\ref{Thm:ComplexityOne} changes  as the tolerance decreases, on the example data $\kappa = 1000$, $f(x_0)-f^* = 10^5$ and $\omega = 1$.}
    \label{fig:ItComp}
\end{figure}

The results of Theorem~\ref{Thm:ComplexityOne} are illustrated in Figure~\ref{fig:ItComp}, via a numerical example with $\kappa = 1000$, $f(x_0)-f^* = 10^5$, and assume that relaxation is never used ($\omega = 1$). Recall that \eqref{eq:K} is independent of $\ell$, which gives the same iteration complexity bound for all Relaxed $\ell$-MGD methods. This is represented by the black `dash-dot' line in Figure~\ref{fig:ItComp}, which shows the iteration complexity bound $K$ as the tolerance $\epsilon$ ranges from $0.1$ to $10^{-8}$. Unsurprisingly, $K$ is larger when $\epsilon$ is smaller, i.e., as the stopping tolerance decreases, the number of iterations increases. On the other hand \eqref{eq:Khat} depends upon $\ell$. The solid lines in Figure~\ref{fig:ItComp} show this complexity bound $K$ as $\epsilon$ varies, with each of the coloured lines corresponding to a specific Relaxed $\ell$-MGD method.

Notice that the slope of the lines correspond to the rate factor, so the slopes for \eqref{eq:Khat} are all the same, and they are steeper/better than the slope for \eqref{eq:K}. This figure also shows that using \eqref{eq:Khat}, the bound $K$ is larger whenever $\ell$ is larger. It can be seen that if $\ell \leq 3$, then \eqref{eq:Khat} should be used to compute the bound $K$, if $\ell \geq 6.5$, then \eqref{eq:K} should be used, while if $3< \ell < 6.5$, the `better' bound will depend upon the choice of the stopping tolerance $\epsilon$.

\section{Computational Practicalities}\label{sec:EfficientComp}

It is widely accepted that algorithms satisfying Definition~\ref{def:Algorithms}, with even moderate values of $\ell$, are computationally intractable, because the step size \eqref{alpha} (which must be computed at every iteration) involves a large matrix power. However, it is now shown that, with an initial overhead of $\lfloor \ell \rfloor+2$ matrix-vector products, and additional storage capacity of $\lfloor \ell \rfloor+2$ vectors, each algorithm in the family can be implemented using only one matrix-vector product per iteration. 

Note that successive gradients can be updated as 
\begin{equation}\label{gradientupdate}
    g\kpo = (I - \alpha\kk A)g\kk,
\end{equation}
while the $y\kk$s can be updated (in exact arithmetic) similarly,
\begin{eqnarray}\label{updateruleyk}
y\kpo 
&\overset{\eqref{y}}{=}& A^{\ell}g\kpo\notag\\
&\overset{\eqref{gradientupdate}}{=}& A^{\ell}(I - \alpha\kk A) g\kk\notag\\
&=& (I - \alpha\kk A) A^{\ell}g\kk\notag\\
&=& (I - \alpha\kk A) y\kk.
\end{eqnarray}

Given an initial point $x\kz\in \R^n$, and corresponding initial gradient $g\kz\in \R^n$, define $v\kz^{(j)} \eqdef A^j g\kz $ for $j = 0,\dots, \lfloor \ell \rfloor + 1$. These $\lfloor \ell \rfloor + 2$ vectors can be initialised efficiently (and then stored) using the recursion:
\begin{eqnarray}\label{vjupdategeneral}
    v\kz^{(j)} = 
    \begin{cases}
        g\kz & \text{if }j = 0\\
        A v\kz^{(j-1)} \;(\equiv A^jg\kz) & \text{for }j = 1,\dots, \lfloor \ell \rfloor + 1.
    \end{cases}
\end{eqnarray}
Then, for all $k\geq 0$, given \eqref{vjupdategeneral}, (i.e., given $v_k^{(0)},\dots,v_k^{(\lfloor \ell \rfloor + 1)}$), $\alpha\kk$ can be computed via
\begin{equation}\label{alphaupdate}
\alpha\kk = \frac{(v\kk^{(\lfloor \ell \rfloor)})^Tv\kk^{(\lfloor \ell +\frac12 \rfloor)}}{(v\kk^{(\lfloor \ell +\frac12 \rfloor)})^Tv\kk^{(\lfloor \ell + 1 \rfloor)}}.
\end{equation}
Furthermore, for all $k\geq 1$, one can update
\begin{eqnarray}\label{vjupdatekpo}
    v\kpo^{(j)} = 
    \begin{cases}
        v\kk^{(j)} - \alpha\kk v\kk^{(j+1)}, & \text{for }j = 0,\dots, \lfloor \ell \rfloor\\
        A v\kk^{(\lfloor \ell \rfloor)}  & \text{if }j = \lfloor \ell \rfloor + 1.\\
    \end{cases}    
\end{eqnarray}
Thus, it is only when $k=0$ that the initial $j=\lfloor \ell \rfloor + 1$ vectors need be computed explicitly, and then for $k\geq 1$  all vectors but one can be updated recursively.

\begin{remark}
Notice that the scheme \eqref{vjupdategeneral} --- \eqref{vjupdatekpo} does not involve any matrix square roots.
\end{remark}

To make the ideas above concrete, and to verify the formulae, consider the following example.
\begin{example}[Implementation example for $\ell = 5/2$]
    For $k=0$, given $x\kz\in \R^n$, recall \eqref{vjupdategeneral} and note that $\lfloor \ell \rfloor  = \lfloor 5/2 \rfloor = 2$. So, compute and store $v\kz^0 = g\kz$, $v\kz^1 = Av\kz^0 = g\kz$, $v\kz^2 = Av\kz^1 = A^2 g\kz$ and $v\kz^3 = Av\kz^2 = A^3 g\kz$. Now
    \begin{equation}
        \alpha\kz \overset{\eqref{alphaupdate}}{=} \frac{(v\kk^2)^Tv\kk^3}{(v\kk^3)^Tv\kk^3}, \overset{\eqref{vjupdategeneral}}{=} \frac{g\kz^TA^5g\kz}{g\kz^TA^6 g\kz} \overset{\eqref{y}}{=} \frac{y\kz^Ty\kz}{y\kz^TA y\kz}. 
    \end{equation}
    For $k=1$ update $x_1 = x\kz - \alpha\kz g\kz$. Now,
    \begin{eqnarray*}
        v_1^0 &=& v_0^0 - \alpha_0 v_0^1 = g_0 - \alpha_0 Ag_0 = g_1,\\
        v_1^1 &=& v_0^1 - \alpha_0 v_0^2 = Ag_0^0 - \alpha_0 A^2g_0 = A(I - \alpha_0 A)g\kz = Ag_1\\
        v_1^2 &=& v_0^2 - \alpha_0 v_0^3 = A^2g_0^0 - \alpha_0 A^3g_0 = A^2(I - \alpha_0 A)g\kz = A^2g_1
    \end{eqnarray*}
    and compute explicitly $v_1^3 = Av_1^2 = A^3g_1 $. Thus, only one matrix vector product is calculated, as expected.
\end{example}

\begin{remark}
It is well known that updating strategies (such as the one described above) should be used with caution, due to the potential for numerical instabilities/inaccuracies to arise (and propagate). In such cases, it may be pertinent to employ restarts to mitigate the risk, e.g., every $q$ iterations, say, the matrix vector products are computed from scratch. A thorough investigation into such issues is left for future research.
\end{remark}

\section{Numerical Experiments}\label{sec:NumericalResults}

In this section, numerical experiments are presented to demonstrate the practical behaviour of several of the algorithms discussed in this work. All the code is written in Python 3.7.4, using an AMD Ryzen 5 3600 CPU with 16GB of RAM.\footnote{Note that in all the figures the y-axis is $\text{log}_{10}$--scale.}
\footnote{The phrase `randomly generated' is used to mean that the entries of the matrix/vector are uniformly distributed random entries from the half open unit interval $[0,1)$.} 

\subsection{Selection of the relaxation parameter}
\label{Sec:numericalomega}
To the best of our knowledge, the work of Raydan and Svaiter in \cite{Raydan2002} is the first to show the effects of relaxation on the Steepest Descent method (i.e., $\ell = 0$). There, at each iteration $k\geq 0$ the relaxation parameter $\omega\kk \in (0, 2)$ is chosen randomly, whereas in this work the relaxation parameter is fixed for all $k$. Two questions arise: (1) `\emph{how should the relaxation parameter $\omega$ be chosen in practice?}', and (2) `\emph{how does random relaxation compare with fixed relaxation?}'. These questions are investigated now, and it seems pertinent to use an identical experimental set up to that in \cite[Section~3]{Raydan2002}.

To this end, $A \in \R^{1000 \times 1000}$ is diagonal, where the $i$th diagonal entry is $A_{ii} = i$, for $i=1,\dots,1000$, giving a condition number of $\kappa(A) = 1000$. Further, $b = 0$, so that $x^* = 0$ and $f^* = 0$. The starting iterate $x_0\in \R^n$ is randomly generated and the stopping condition is $\|g\kk\|_2^2 < 10^{-8}$ (or 1000 iterations). 

Three algorithms are applied to this problem: Relaxed $\ell$-MGD methods (recall Definition~\ref{def:Algorithms}) where $\ell \in \{0,1/2,1\}$. All algorithms were run multiple times using a range of different (fixed) relaxation values $\omega \in (0,2)$, and all trials were repeated multiple times using different initial vectors $x\kz$.

\begin{figure}[ht!]
    \centering
    \includegraphics[width =0.52\textwidth]{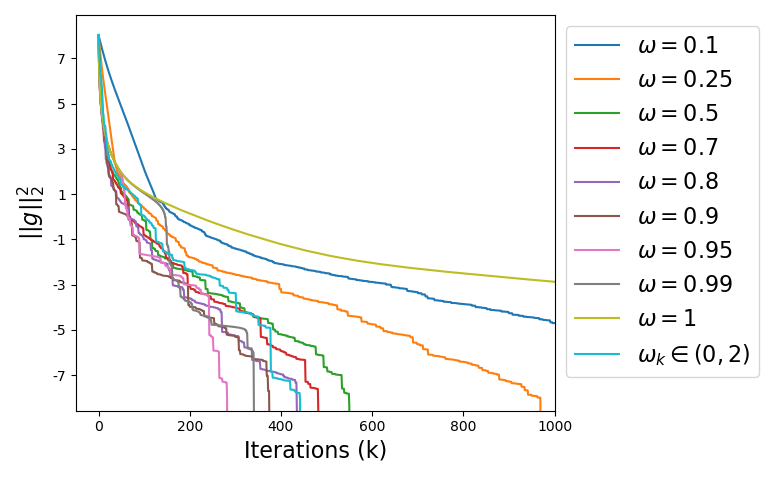}
    \includegraphics[width =0.42\textwidth]{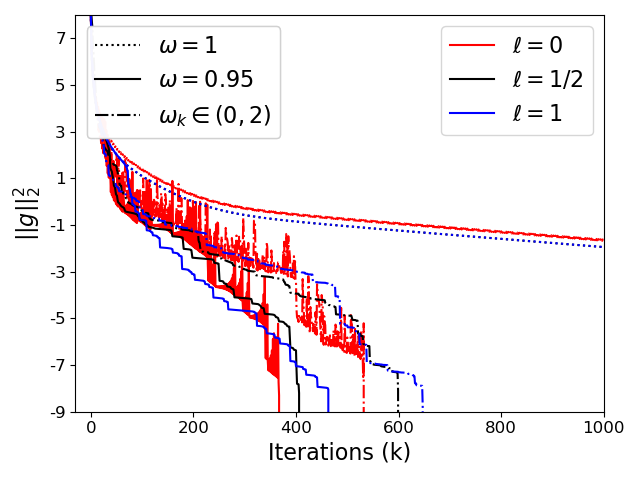}
    \caption{Left plot showing the evolution of the gradient for various values of the relaxation parameter $\omega$ on the experiment described in Section~\ref{Sec:numericalomega}. Right plot showing the evolution of the (squared) gradient norm as the Relaxed $\ell$-MGD method (used with either $\ell=0,1/2,1$, and specified relaxation) is run on the same experiment.}
    \label{fig:omegacomparison}
\end{figure}

Figure~\ref{fig:omegacomparison} (left plot) shows the results of one such trial, using the Relaxed $1$-MGD method, with $\omega\in \{0.1,0.25,0.5,0.7,0.8,0.9,0.95,0.99,1\}$. This instance was chosen because it exhibits behaviour that was representative of that observed over all trials and all algorithms ($\ell \in \{0,1/2,1\}$). The case $\omega = 1$ (no relaxation) is included as a control, to allow comparison with relaxation. Notice that all shown choices of the relaxation parameter resulted in a reduction in the number of iterations compared with no relaxation. Further, when $0.5\leq \omega <1$, all runs achieved the stopping condition in fewer than 600 iterations. This demonstrates a clear improvement in performance over a wide range of relaxation parameter values compared with the baseline case. Also shown in Figure~\ref{fig:omegacomparison} (left plot) is an instance of random relaxation (i.e., the Relaxed $1$-MGD method with $\omega\kk \in (0,2)$ for all $k\geq 0$). This shows the very good practical performance of random relaxation, but notice that the runs with fixed $\omega\in \{0.9,0.95,0.99\}$ were all better (fewer iterations).

Experiments with relaxation values larger than 1 (i.e., $\omega\in (1,2)$) are omitted, as they were observed to be slower (more iterations) than the control case $\omega = 1$. This agrees with the findings in van den Doel and Ascher \cite[Section~5]{Doel2012}.

These results demonstrate that (i) relaxation can help the practical performance of the methods in Definition~\ref{def:Algorithms}; and that (ii) the improvement in performance is not too sensitive to the specific choice of $\omega$ value (fixed relaxation values in the range $\omega \in [0.9,0.99]$ resulted in strong performance). As the choice $\omega =0.95$ worked consistently well, this value is adopted for the remaining numerical experiments.

The right plot in Figure~\ref{fig:omegacomparison} shows the results of another trial run for this experiment. Here, the evolution of gradient norm is shown for all three algorithms ($\ell = 0,1/2,1$), and for three values of the relaxation parameter: (i) no relaxation $\omega = 1$, (ii) constant relaxation $\omega = 0.95$, and (iii) random relaxation $\omega\kk \in (0, 2)$ for $k\geq 0$. It can be seen that if relaxation is used (either fixed or random), all algorithms met the stopping condition in fewer than 800 iterations. In contrast, when $\omega = 1$, all algorithms were terminated after reaching the maximum number of allowed iterations, at which point the gradient 2-norm was still relatively large ($\|g_{1000}\|_2^2 \approx 0.1$). Again, this clearly highlights the benefits of including relaxation when using minimal gradient type methods.

The right plot in Figure~\ref{fig:omegacomparison} also provides insight into the difference between random versus fixed relaxation. When random relaxation is used ($\omega_k \in (0,2)$) all algorithms met the stopping condition in around 700 iterations, while if a fixed relaxation parameter ($\omega=0.95$) is used then approximately 400 iterations are needed. This shows that a well-chosen fixed relaxation parameter can be a good choice and outperform random relaxation.

\subsection{Experiments using synthetic datasets}\label{sec:RandommatrixTest}

In this section the algorithms (recall Definition~\ref{def:Algorithms}) are applied to \eqref{f}, where the data is randomly generated (synthetic). The matrix $A\in \R^{n\times n}$ is either dense or sparse (see details below), $x_0\in \R^n$, and $x^*\in \R^n$ are randomly generated, $b = A\xs$ and the stopping condition is $\|g\kk\|_2^2\leq 10^{-6}$.

\paragraph{Dense random matrices.} Here, the positive definite matrix $A$ is generated as follows. Let $m=1500$, let $n=1000$, let $B \in \R^{m \times n}$ be randomly generated and form $A = B^TB$. Observe that $A$ is symmetric, and because $m>n$ and $B$ is dense, it is highly likely that $A$ is positive definite (i.e., Assumption~\ref{Assume1} is satisfied).

\begin{figure}
    \centering
    \includegraphics[height=0.3\textheight]{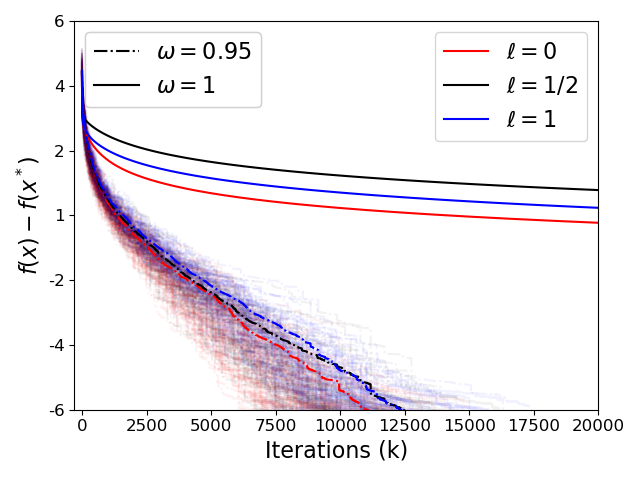}
    \caption{The results of the experiment described in Section~\ref{sec:RandommatrixTest} in the dense case with $A \in \R^{n \times n}$ where $n = 10^3$. The relaxed $\ell$-MGD algorithms for $\ell \in \{0, 1/2, 1\}$ with $\omega = 1$ and $\omega = 0.95$ are tested. The evolution of the function values is displayed. The mean over all runs for each of the methods are shown in bold, and the pale lines correspond to the trajectories for each of the 100 individual instances for each method.}
    \label{fig:DenseTests}
\end{figure}

The Relaxed $\ell$-MGD methods with $\ell = 0, 1/2, 1$ are applied to this problem, both with relaxation ($\omega = 0.95$) and without ($\omega = 1$), and the experiment was repeated 100 times. Figure~\ref{fig:DenseTests} shows the evolution of the function value difference for this experiment. There are a total of 300 pale `dash-dot-line' curves in this plot; 3 curves for each of the $\ell = 0,1/2,1$ algorithms, and for each of the 100 trial runs. Overlaid on these results are 3 darker/bold lines, corresponding to the average function value evolution over all 100 runs for each algorithm/$\ell$-value. Also shown in the figure is 3 solid lines corresponding to the average function value evolution over all 100 runs for each algorithm/$\ell$-value without relaxation (the individual runs are not shown because no run was competitive with relaxation).

It can be seen in Figure~\ref{fig:DenseTests} that when $\omega=1$, on average all three algorithms perform similarly, with Steepest Descent ($\ell=0$) being slightly better than the other methods. It can also be seen that there is a large improvement in performance when relaxation is used, with all algorithms requiring far fewer iterations, and on average the termination criteria was met within 12500 iterations. Again, all algorithms perform similarly on average when $\omega = 0.95$, with Steepest Descent ($\ell=0$) being slightly better.\footnote{Note that a plot showing the evolution of $\|g\kk\|_2^2$ is omitted. This is because the oscillations in the gradient norm in the $\ell=0$ case overwhelm the remaining curves making it difficult to `see' the results. }

\paragraph{Sparse random matrices.}
Here $A\in \R^{n\times n}$ ($n=10^6$) is sparse (the scipy.sparse package was used) and to ensure that Assumption~\ref{Assume1} is met, $A$ is constructed as follows. Let $C\in \R^{n\times n}$ be a sparse matrix with an average of 10 nonzero randomly generated entries per row and let $B = C + C^T$  so that $B$ is symmetric. Now, let $z\in \R^n$ be randomly generated with uniform entries on $(0,1000)$, and let $e$ denote the vector of all ones. Then $A = B + {\rm diag}(Be+z)$. Constructing $A$ in this way ensures positive-definiteness with strict diagonal dominance. As in previous experiments, the Relaxed $\ell$-MGD algorithms with $\ell = 0,1/2,1$ are applied to this problem set up, both with relaxation ($\omega = 0.95$) and without ($\omega = 1$), and the experiment was repeated 10 times.

\begin{figure}[ht!]
    \centering
    \includegraphics[width=0.49\textwidth]{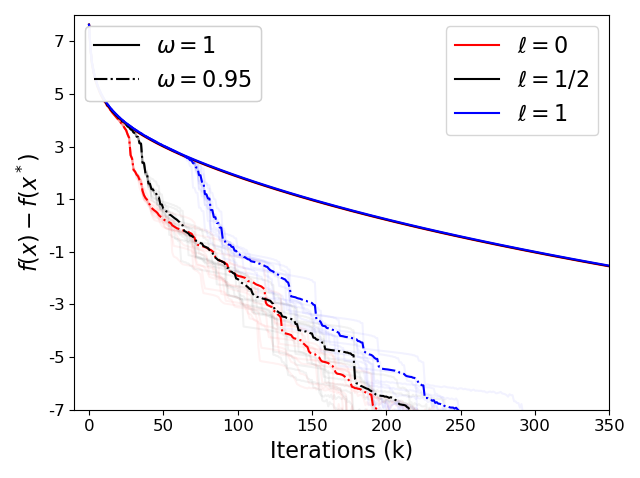}
    \includegraphics[width=0.49\textwidth]{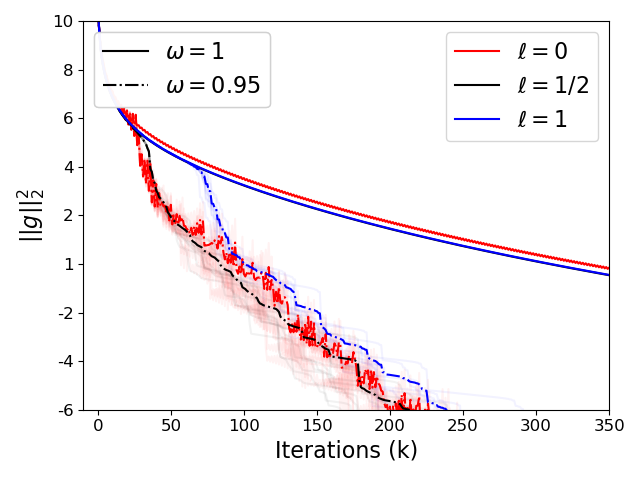}
    \caption{The results of the experiment described in Section~\ref{sec:RandommatrixTest} in the sparse matrix case with $A \in \R^{n \times n}$ where $n = 10^6$. The relaxed $\ell$-MGD algorithms for $\ell \in \{0, 1/2, 1\}$ with $\omega = 1$ and $\omega = 0.95$ are tested. The evolution of the function values is shown in the left plot and the evolution of the 2-norm squared of the gradient is shown on the right. The mean over all runs for each of the methods are shown in bold, and trajectories for each of the 10 instances of each method are shown in pale.}
    \label{fig:Sparse1mill}
\end{figure}

The results are shown in Figure~\ref{fig:Sparse1mill} and the curves are similar to the dense case, with the pale lines corresponding to individual runs, while the bold lines correspond to averages over all trial runs. The left plot shows the evolution of the function values while the right plot shows the norm of the gradient.\footnote{Fewer trials are shown here (10 versus 100 in the dense case), so that the evolution of the gradient could be included without the oscillations for the $\ell=0$ case masking the other trials. However, experiments using 100 runs were also performed, and the results (omitted) are similar to that shown in Figure~\ref{fig:Sparse1mill}.} Again, it can be seen that algorithm performance improves markedly when relaxation is used; without relaxation approximately 1000 iterations were needed to achieve the desired stopping condition, whereas the algorithms terminated in 200--300 iterations with relaxation. 
Figure~\ref{fig:Sparse1mill} shows the evolution of the function values and norm of the gradient as the algorithms progress.

\subsection{Experiments using LIBSVM datasets}\label{sec:LIBSVMDatasetTests}

In this section three datasets from LIBSVM \cite{LIBSVMCC01a} are considered: \texttt{mushrooms}, \texttt{a1a} and \texttt{w1a}. Basic information about the datasets is given in Table~\ref{tab:datasetsingularvals}, where $\sigma_1$ and $\sigma_n$ denote the largest and smallest singular values of the data matrix, respectively.
\begin{table}[ht!]
    \centering
    \begin{tabular}{||c|c|c|c|c||}
    \hline
        Dataset & Rows & Columns & $\sigma_n$ & $\sigma_1$ \\
        \hline
        mushrooms & 8124 & 112 & 1e-14 & 289\\
        a1a & 1605 & 119 & 1e-15 & 100.3\\
        w1a & 2477 & 300 & 2e-15 & 78.5\\\hline
    \end{tabular}
    \caption{Details for the \texttt{mushrooms}, \texttt{a1a} and \texttt{w1a} LIBSVM datasets.}
    \label{tab:datasetsingularvals}
\end{table}
Note that for each of the datasets, the data matrix --- denoted by $B$ --- is rectangular, and let $d$ denote the corresponding labels vector. The following regularized problem is considered here:
\begin{equation}\label{regprob}
    \min_x \|Bx-d\|_2^2 + \lambda \|x\|_2^2,
\end{equation}
where the regularization parameter is chosen to be $\lambda = 10^{-6}$. Now, (ignoring the constant), \eqref{regprob} fits \eqref{f} with $A \gets B^TB + \lambda I$ (which satisfies Assumption~\ref{Assume1}) and $b \gets B^Td$. The singular values stated in Table~\ref{tab:datasetsingularvals} show that matrices, $A$, considered here are ill-conditioned, despite the regularization. 

\begin{figure}[ht!]
    \centering
    \includegraphics[width=0.49\textwidth]{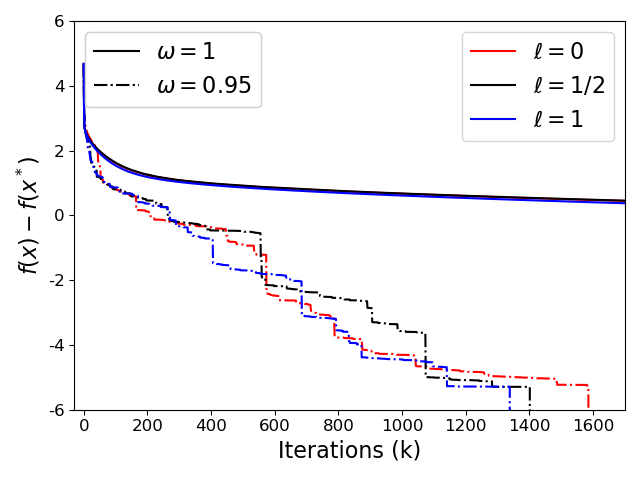}
    \includegraphics[width=0.49\textwidth]{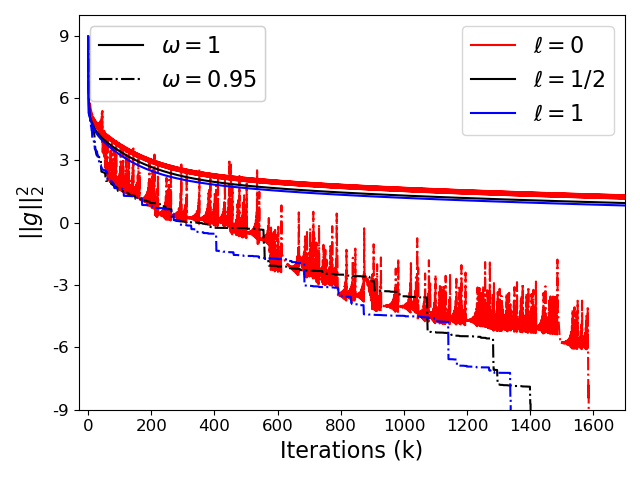}
    
    \includegraphics[width=0.49\textwidth]{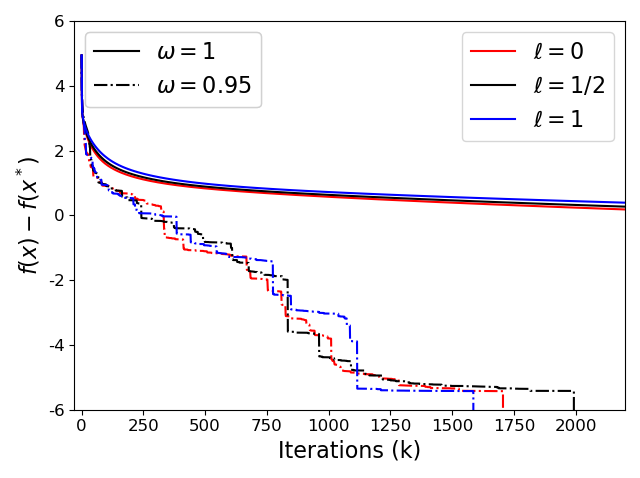}
    \includegraphics[width=0.49\textwidth]{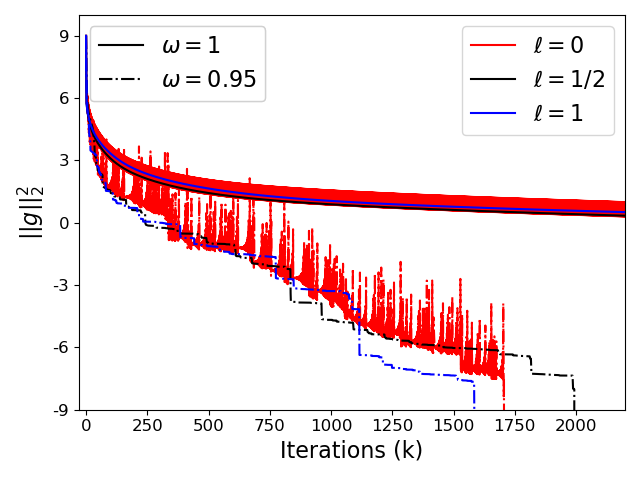}
    
    \includegraphics[width=0.49\textwidth]{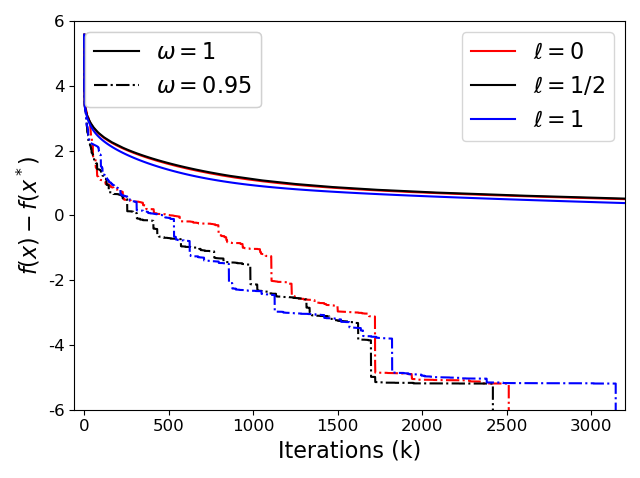}
    \includegraphics[width=0.49\textwidth]{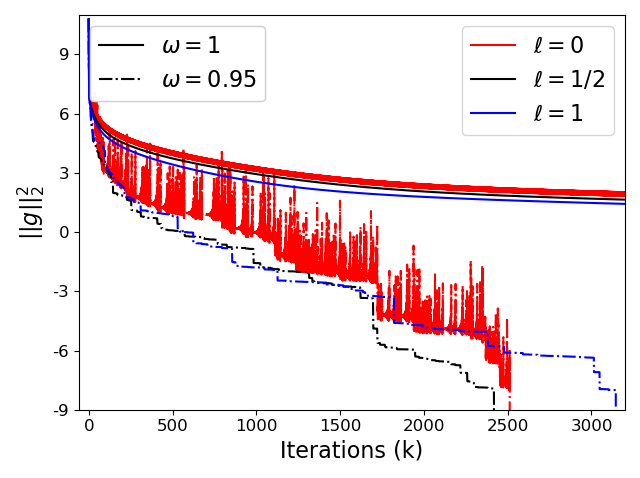}
    \caption{Results for the numerical experiments described in Section~\ref{sec:LIBSVMDatasetTests}. Rows 1--3 correspond to the datasets  \texttt{a1a}, \texttt{w1a}, and \texttt{mushrooms}, respectively. The left column shows the evolution of the function values for each of the algorithms and the right column shows the evolution of the 2-norm of the gradient as iterates progress.}
    \label{fig:datasetplots}
\end{figure}

Three Relaxed $\ell$-MGD methods ($\ell \in \{0, 1/2, 1\}$) were employed on the experimental set up just described, both with ($\omega = 0.95$) and without ($\omega = 1$) relaxation. Note that pythons in-built least squares solver was used to compute $x^*$, and the stopping condition is $\|g\kk\|_2^2<10^{-9}$. Multiple trials were run, where in each trial all algorithms were initialized using the same randomly generated initial point $x_0$, and a representative example of the behaviour of these methods is shown in Figure~\ref{fig:datasetplots}.

Figure~\ref{fig:datasetplots} shows that for each dataset, all algorithms (for $\ell \in \{0, 1/2, 1\}$) perform similarly, and that is typical of the behaviour observed over the many trial runs. Notice that when relaxation is used, the algorithms all terminated in around 2000 iterations because they met the stopping conditions, whereas when there is no relaxation ($\omega = 1$), the error in the function values or gradient norm is several orders of magnitude bigger.

\section{Conclusion}

This work presented new theoretical results for the family of Relaxed $\ell$-Minimal Gradient Descent methods for quadratic optimization. By choosing the norm appropriately, results can be stated concisely for all members of the family, including that (i) the norm of the gradient; (ii) the distance of the iterates from optimality; and (iii) the function values, all converge linearly. A counterexample showed that it is not possible to extend the function value convergence rate for SD \eqref{ratenorelax} to the remaining members of the family. The results are nonasymptotic, so it was also possible to establish iteration complexity results for the algorithms. Moreover, all theoretical results hold when fixed relaxation is employed.

It was also shown that, given a fixed overhead and storage budget, every Relaxed $\ell$-Minimal Gradient Descent method can be implemented using exactly one matrix vector product. This shows that, Relaxed $\ell$-Minimal Gradient Descent methods with large $\ell$ values are not intractable, as was previously thought.

Finally, numerical experiments were presented, which demonstrate the benefits of relaxation for this family of methods. The results were consistent with those found in the existing literature, and support the widely held view that a fixed relaxation value of $\omega = 0.95$ works well in practice.

\paragraph{Future work.} This work describes the non-asymptotic behaviour for the Relaxed $\ell$-Minimal Gradient Descent methods, although a loose end remains. In particular, the counterexample in Section~\ref{sec:Counterexample} shows that it is not possible to extend the convergence rate for SD to the rest of the family. However, the authors have not observed the `two step' decrease in the function value to be worse than the Kantorovich rate numerically, and it is an open question as to whether it is possible to prove a `two step' ($f_{k+2}-f\kk$ vs $f\kpo-f\kk$) type result (see also \cite[Section~5]{Nocedal06}).

\bibliographystyle{unsrt}
\bibliography{refs}

\appendix

\section{Comparison with the results in \cite{Pronzato_2005}}\label{app:PronzatoEquiv}
The rate of convergence for the norm of the gradient stated in Theorem~\ref{Thm:gmonodecreaserelaxed} is equivalent to the rate of convergence given by Pronzato et al. in \cite{Pronzato_2005}. That work considers a class of methods, called the $P$-gradient algorithms, that are used to solve problem \eqref{f}. These methods are defined as follows.
\begin{definition}[Definition 1 in \cite{Pronzato_2005}]
   Let $P(\cdot)$ be a real function defined on $[m, M]$, infinitely differentiable, with Laurent series
     $P(z) = \sum_{-\infty}^{\infty} c\kk z^k, c\kk \in \R \text{ for all } k,$
   such that $0 < \sum_{-\infty}^{\infty} c\kk a^k < \infty $ for $a \in [m, M]$. The $k$-th iteration of a $P$-gradient algorithm is defined by 
   \begin{equation}
       x\kpo = x\kk -\gamma_k g\kk
   \end{equation}
   where the step-length $\gamma\kk$ minimises $\|g\kpo\|_{P(A)}^2$ with respect to $\gamma$, with $g\kpo =  \nabla f (x\kk - \gamma g\kk)$.
\end{definition}

The work \cite{Pronzato_2005} uses the setting of a Hilbert space, whereas this work considers $\R^n$. Thus, here the following notation translation can be made:
\begin{eqnarray*}
    P(A) = A^{2\ell-1},\qquad m = \lambda_n,\qquad \text{and}\qquad M = \lambda_1.
\end{eqnarray*}
Considering \cite[Equation (9)]{Pronzato_2005} used with the notation of this work shows that
\begin{equation}\label{Prozantorate}
    r\kk = \frac{\|g\kpo\|^2_{A^{2\ell-1}}}{\|g\kk\|^2_{A^{2\ell-1}}} = 1 - \frac{1}{L_k}.
\end{equation}
Now, \cite[(5) and (6)]{Pronzato_2005} become
\begin{equation}\label{zandmu}
    z_k \eqdef \frac{(A^{2\ell-1}\cdot A)^{1/2}g\kk}{\| g\kk \|_{(A^{2\ell-1}\cdot A)}} \equiv \frac{A^{\ell}g\kk}{\| A^{\ell}g\kk \|_{2}} \overset{\eqref{y}}{=} \frac{y\kk}{\| y\kk \|_{2}}\quad \text{and}\quad \mu^k_j \eqdef \|z\kk \|_{A^{j}}^2, \;\; j\in \mathbb{Z}.
\end{equation}
Finally, the expression for $L_k$ is given in \cite[p.414]{Pronzato_2005}, so that combining with \eqref{zandmu} gives
\begin{eqnarray}
    L_k & =& \mu_1^k \mu_{-1}^k \notag \\
    &=& \|z\kk \|^2_{A} \|z\kk \|^2_{A^{-1}} \notag\\
    &=& \frac{ \| y\kk \|^2_A}{\|y\kk\|^2_2} \frac{\| y\kk \|^2_{A^{-1}} }{\|y\kk\|^2_2} \notag\\
    &=& \frac{\|y\kk\|_{A}^2\|y\kk\|_{A^{-1}}^2}{\|y\kk\|^4_2}.\label{Lk}
\end{eqnarray}
It remains to note that substituting \eqref{Lk} into the rate expression \eqref{Prozantorate} gives \eqref{PronzatoEquiv} in the proof of Theorem~\ref{Thm:gmonodecreaserelaxed} (when $\omega = 1$), i.e., no relaxation).

\section{Technical Lemmas for Section~\ref{Sec:functionvals}}
\label{sec:technicallemmas}

Here, several technical lemmas are established, which are required for the proof of Theorem~\ref{thm:fvaldecrease}.

\begin{lemma}[Fact 8.12.7 in \cite{Bernstein05}]
Let $A\in \R^{n\times n}$ be positive semidefinite and let $u\in \R^n$. Then $(u^TA^2u)^2 \leq (u^TAu)(u^TA^3u)$
and
\begin{equation}\label{fact2}
    (u^TAu)^2 \leq (u^Tu)(u^TA^2u).
\end{equation}
\end{lemma}

\begin{lemma}\label{Lemmainduction}
    Let $f$ be given in \eqref{f}, let Assumption \ref{Assume1} hold, and let $\ell \in \{0,\tfrac12,1,\tfrac32,2,\tfrac52,3,\dots\}$ be fixed. Given a point $x\kk\in\R^n$, let $g\kk$ and $y\kk$ be defined in \eqref{gradient} and \eqref{y}, respectively. Then 
    \begin{equation}\label{intermediateThmresult}
        \frac{\|y\kk\|_2^2}{\|y\kk\|_A^2}\frac{\|g\kk\|_A^2}{\|g\kk\|_2^2} \leq 1.
    \end{equation}
\end{lemma}
\begin{proof}
    Note that \eqref{intermediateThmresult} is equivalent to 
    \begin{equation}\label{inductionhypothesisint}
        \frac{g\kk^TA^{2\ell}g\kk}{g\kk^TA^{2\ell+1}g\kk}\frac{g\kk^TAg\kk}{g\kk^Tg\kk} \leq 1,
    \end{equation}
    which holds with equality when $\ell = 0$. Now, the induction hypothesis is to assume that \eqref{inductionhypothesisint} holds for some $\ell \geq 1/2$.
    Substituting $u = A^{\ell}g\kk$ into \eqref{fact2} shows that
    \[(g\kk^TA^{2\ell+1}g\kk)^2 \leq (g\kk^TA^{2\ell}g\kk)(g\kk^TA^{2\ell+2}g\kk).
    \]
    Rearranging and then multiplying through by $g\kk^TAg\kk/g\kk^Tg\kk$ gives
    \begin{eqnarray}
        \frac{g\kk^TA^{2\ell+1}g\kk}{g\kk^TA^{2\ell+2}g\kk}\frac{g\kk^TAg\kk}{g\kk^Tg\kk} \leq\frac{g\kk^TA^{2\ell}g\kk}{g\kk^TA^{2\ell+1}g\kk}\frac{g\kk^TAg\kk}{g\kk^Tg\kk} \overset{\eqref{inductionhypothesisint}}{\leq} 1.
    \end{eqnarray}
    Hence, \eqref{intermediateThmresult} is true for all $\ell \in \{0,\tfrac12,1,\tfrac32,2,\tfrac52,3,\dots\}$, which completes the proof.
\end{proof}

\begin{lemma}\label{lemmaBounding}
    Let $A\in \R^{n\times n}$ and let Assumption~\ref{Assume1} hold. Let $\ell \in \{0,\tfrac12,1,\tfrac32,2,\tfrac52,3,\dots\}$ be fixed. Then, for any nonzero vector $z\in \R^n$, the following inequalities hold,
    \begin{eqnarray}\label{ineq1}
  z^TAz \leq  \frac{1}{\lambda_n^{2\ell}}\;  \|z\|_{A^{2\ell+1}}^2, 
\end{eqnarray}
and
\begin{eqnarray}\label{ineq2}
  \|z\|_{A^{2\ell+1}}^2 \leq  \lambda_1^{2\ell} \; z^TAz.
\end{eqnarray}
\end{lemma}
\begin{proof}
   By the Spectral Theorem, every symmetric matrix is orthogonally diagonalizable. So, if $A$ has eigenvalues $0< \lambda_n(A)\leq \cdots\leq \lambda_1(A)$, then $A^{2\ell}$ has eigenvalues $\lambda_i(A^{2\ell}) = (\lambda_i(A))^{2\ell}$ for $i=1,\dots,n$. Now, taking $B = A^{2\ell}$, and $v = A^{1/2}z$ in \eqref{RayleighQuotient} shows that
\begin{eqnarray*}
  \lambda_n^{2\ell} \; z^TAz \leq  z^TA^{2\ell+1}z  \leq   \lambda_1^{2\ell} \; z^TAz.
\end{eqnarray*}
Noting $z^TA^{2\ell+1}z \;\equiv \|z\|_{A^{2\ell+1}}^2$ gives \eqref{ineq1} and \eqref{ineq2}.
\end{proof}

\end{document}